\documentclass{article}

\usepackage{arxiv}

\usepackage[utf8]{inputenc} 
\usepackage[T1]{fontenc}    
\usepackage{hyperref}       
\usepackage{url}            
\usepackage{booktabs}       
\usepackage{amsfonts}       
\usepackage{nicefrac}       
\usepackage{microtype}      
\usepackage{lipsum}		
\usepackage{graphicx}
\usepackage{natbib}
\usepackage{doi}
\usepackage{amsmath, amsthm, amssymb}
\newcommand{\bra}[1]{\left(#1\right)}
\newcommand{\sbra}[1]{\left[#1\right]}
\newcommand{\seq}[1]{\left<#1\right>}
\newcommand{\vertiii}[1]{{\left\vert\kern-0.25ex\left\vert\kern-0.25ex\left\vert #1
    \right\vert\kern-0.25ex\right\vert\kern-0.25ex\right\vert}}

\newcommand{\abs}[1]{\left\vert#1\right\vert}
\newcommand{\set}[1]{\left\{#1\right\}}
\renewcommand{\c}{\mathbb C}
\newcommand{\Real}{\mathbb R}

\renewcommand{\to}{\longrightarrow}
\newcommand {\LEF}{\mathcal{L}(\mathfrak{E},\mathfrak{F})}

\newcommand {\LE}{\mathcal{L}(\mathfrak{E})}
\renewcommand {\L}{\mathcal{L}}

\newcommand{\A}{\mathfrak{A}}
\newcommand{\E}{\mathfrak{E}}
\newcommand{\F}{\mathfrak{F}}

\renewcommand{\t}{\mathbf{t}}
\newcommand{\q}{\mathbf{q}}
\newcommand{\x}{\mathbf{x}}

\newcommand{\y}{\mathbf{y}}

\newtheorem{theorem}{Theorem}[section]
\newtheorem{lemma}[theorem]{Lemma}

\newtheorem{corollary}[theorem]{Corollary}
\newtheorem{definition}[theorem]{Definition}

\newtheorem{remark}[theorem]{Remark}

\newcommand\mystyle{\everymath{\displaystyle}}
\mystyle
\title{Some inequalities for the Euclidean operator radius
of two operators in Hilbert $C^{\ast}$-Modules  space}

\author{ \href{https://orcid.org/0000-0002-3816-5287}{\includegraphics[scale=0.06]{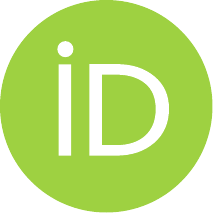}\hspace{1mm}M.H.M.~Rashid}\thanks{Use footnote for providing further
		information about author (webpage, alternative
		address)---\emph{not} for acknowledging funding agencies.} \\
	Department of Mathematics\&Statistics\\Faculty of Science P.O.Box(7)\\
	Mu'tah University University\\
	Mu'tah-Jordan \\
	\texttt{mrash@mutah.edu.jo}
}

\renewcommand{\shorttitle}{\textit{arXiv} Template}

\hypersetup{
pdftitle={Some inequalities for the Euclidean operator},
pdfsubject={q-bio.NC, q-bio.QM},
pdfauthor={M.H.M.Rashid},
pdfkeywords={numerical radius, inner product space, $C^*$-algebra, $A$-module},
}

\begin{document}
\maketitle

\begin{abstract}
	The Euclidean operator radius of two bounded linear operators in the Hilbert $C^*$-module over $\A$ is given some precise bounds. Their relationship to recent findings in the literature that offer precise upper and lower bounds on the numerical radius of linear operators is also established.
\end{abstract}

\keywords{numerical radius\and inner product space\and $C^*$-algebra\and $\A$-module }

\section{Introduction}
A common tool in operator and operator algebra theory is the Hilbert $C^*$-module. They provide a significant class of instances in the theory of operator $C^*$-modules. In addition, the theory of Hilbert $C^*$-modules is a fascinating subject in and of itself. interacting with operator algebra theory and incorporating concepts. It develops from non-commutative geometry and generates findings as well as new issues that draw attention.\\
\indent An extension of the idea of a Hilbert space is the idea of a Hilbert $C^*$-module. Kaplansky \cite{Kap} used such things for the first time. The induced representations of $C^*$-algebras by \cite{Rief} and Paschke's PhD dissertation \cite{Pasc} marked the beginning of the study of Hilbert $C^*$-modules. The theory of operator algebras, operator $K$-theory, group representation theory, and theory of operator spaces all benefit from the employment of Hilbert $C^*$-modules. Additionally, it is utilized to research the Morita equivalence of $C^*$-algebras, the $C^*$-algebra quantum group, and the $C^*$-algebra $K$-theory (\cite{Lance};  \cite{Wegge-Olsen}).\\
\indent A pre-Hilbert module over a $C^*$-algebra $\A$ is a complex linear space $\E$ which is a right $\A$-module (and $\lambda(ax)=(\lambda a)x=a(\lambda x)$, where $\lambda \in \c$, $a \in \A$, and $x \in \E$) equipped with an $\A$-valued inner product
$\seq{\cdot,\cdot}:\E\times\E\to \A$ satisfying:
\begin{enumerate}
  \item [(i)] $\seq{x,x}\geq 0$ and $\seq{x,x}=0$ if and only if $x=0$;
  \item [(ii)] $\seq{x,\alpha y+\beta z}=\seq{x,y}\alpha+\seq{x,z}\beta$;
  \item [(iii)] $\seq{x,ya}=\seq{x,y}a$;
\item [(iv)] $\seq{x,y}=\seq{y,x}^*$,
\end{enumerate}
for all $x,y,z\in\E$, $a\in\A$ and $\alpha,\beta\in\c$.\\
\indent A pre-Hilbert $\A$-module is called a Hilbert $\A$-module or Hilbert $C^*$-module over $\A$, if it is complete with respect to the norm $\vertiii{x}=\vertiii{\seq{x,x}}^{\frac{1}{2}}$.  Suppose that $\E, \F$ are Hilbert $C^*$-modules. We define $\LEF$
to be the set of all maps $t: \E\to \F$ for which there is a map $t^*:\F\to\E$ which satisfies $\seq{tx,y}=\seq{x,t^*y}$ for all $x\in\E$ and $y\in\F$. $\L(\E,\E)$ is simply denoted by $\LE$. It is known that $\LE$ is a $C^*$-algebra. A state on a $C^*$-algebra $\A$ is a positive linear
functional on $\A$ of norm one. We denote the state space of $\A$ by $\varpi(\A)$ (see \cite{Mog}).\\
\indent The major goal of this study is to examine other specific situations of interest and to extend several inequalities in the literature to the Euclidean radius of two operators in a Hilbert $\A$-module. Also offered are conclusions that link the operator norm, the typical numerical radius of a composite operator, and the Euclidean operator radius.

\section{Definitions and Complementary Results}
\begin{definition} Let $\t=(\x,\y)\in\L^2(\E)$. The Euclidean operator radius is
defined by
\begin{eqnarray}\label{Eucw1}
  w_e(\t):&=&w_e(\x,\y) \\
   &=& \sup\set{\bra{\abs{\psi\bra{\seq{\xi,\x \xi}}}^2+\abs{\psi\bra{\seq{\xi,\y \xi}}}^2}^{\frac{1}{2}}:\xi\in\E,\psi\in\varpi(\A)\,\,\mbox{and}\,\,\psi\bra{|\xi|}=1}.\nonumber
\end{eqnarray}
  And the Euclidean operator norm is defined by
  \begin{eqnarray}\label{EucN1}
    \vertiii{\t}:&=&\sqrt{\vertiii{\x^*\x+\y^*\y}} \\
     &=&\sup\set{\bra{\psi\bra{|\x\xi|^2}+\psi\bra{|\y\xi|^2}}^{\frac{1}{2}}:\xi\in\E,\psi\in\varpi(\A)\,\,\mbox{and}\,\,\psi\bra{|\xi|}=1}.\nonumber
  \end{eqnarray}
\end{definition}
\begin{lemma} $\vertiii{\cdot}$ is a norm on $\E$.
\end{lemma}
\begin{proof} If $\t=(t_1,t_2)=0$, then obviously $\vertiii{\t}=0$. If $\vertiii{\t}=0$, then for every $\psi\in\varpi(\A)$ and each $\xi\in\E$ such that $\psi\bra{\abs{\xi}}=1$, we have  $\psi\bra{|t_j\xi|^2}=0$ for all $j=1,2$. We want to show that $\t \xi=0$ for each $x \in \E$.
Fix $\xi\in\E$.
\begin{enumerate}
  \item  If $\psi\bra{|\xi|}=0$, then by the Cauchy-Schwartz inequality we have
$$\psi\bra{\seq{t_j\xi,t_j\xi}}=\psi\bra{\seq{t_j^*t_j\xi,\xi}}\leq \psi\bra{\seq{t_j^*t_j\xi,t_j^*t_j\xi}}^{\frac{1}{2}}\psi\bra{|\xi|}=0,$$
and so $t_j=0$ for all $j=1,2$, i.e., $\t=0$.
  \item If $\psi\bra{|\xi|}\neq 0$, then by taking $y=\frac{x}{\psi\bra{|\xi|}}$, then $\psi\bra{|\zeta|}=1$. By Definition \ref{EucN1},
$\psi\bra{\abs{t_j\zeta}}=0$ for all $j=1,2$ and so $\frac{1}{\phi\bra{|\xi|}}\psi\bra{\abs{t_j\xi}}=0$ for all $j=1,2$.
Thus $\psi\bra{\abs{t_j\xi}}=0$ for all $j=1,2$. Since for every $\psi\in\varpi(\A)$,
we have $\psi\bra{|t_j\xi|}=0$ for all $j=1,2$. We conclude that $|t_j\xi|=0$ for each $\xi \in \E$ and for all $j=1,2$ . So $t_j=0$
 for all $j=1,2$ and so $\t=0$.
\end{enumerate}
\indent On the other hand $\A$ is an abelian $C^*$-algebra, then by \cite[Theorem 3.6]{Jiang}, $|\xi+\zeta|\leq |\xi|+|\zeta|$
for each $\xi,\zeta\in\E$. Thus $\abs{t\xi+s\xi}\leq |t\xi|+|s\xi|$ for all  $\xi\in\E$. Hence
$$\psi\bra{\abs{t\xi+s\xi}}\leq \psi\bra{|t\xi|}+\psi\bra{|s\xi|}.$$
Now by Minkowski's inequality, we have
\begin{eqnarray*}
 \bra{\bra{\psi\abs{(t_1+q_1)\xi}}^2+\bra{\psi\abs{(t_2+q_2)\xi}}^2}^{\frac{1}{2}}&\leq&
\bra{\bra{\psi\bra{|t_1\xi|}+\psi\bra{|q_1\xi|}}^2+\bra{\psi\bra{|t_2\xi|}+\psi\bra{|q_2\xi|}}^2}^{\frac{1}{2}} \\
   &\leq& \bra{\psi\bra{|t_1\xi|^2}+\psi\bra{|t_2\xi|^2}}^{\frac{1}{2}}+\bra{\psi\bra{|q_1\xi|^2}+\psi\bra{|q_2\xi|^2}}^{\frac{1}{2}}
\end{eqnarray*}
Taking the supremum over all $\xi\in\E$ with $\psi\bra{\abs{\xi}}=1$ and $\psi\in\varpi(\A)$, we have
$$\vertiii{\t+\q}\leq \vertiii{\t}+\vertiii{\q}.$$
Finally, if $\t=(t_1,t_2)$, then we have
$$\vertiii{c\t}=\bra{\psi\bra{|ct_1\xi|^2}+\psi\bra{|ct_2\xi|^2}}^{\frac{1}{2}}
=|c|\bra{\psi\bra{|t_1\xi|^2}+\psi\bra{|t_2\xi|^2}}^{\frac{1}{2}}=|c|\vertiii{\t}$$
for all $c\in\c$.
\end{proof}
\begin{theorem}\label{extension} Let $\t=\bra{t_1,t_2}\in\L^n\bra{\E}$.  If $\E$ is a Hilbert $\A$-Modules, then
$$\vertiii{\t}=\sup\set{\bra{\abs{\psi\bra{\seq{\zeta,t_1\xi}}}^2+\abs{\psi\bra{\seq{\zeta,t_2\xi}}}^2}^{\frac{1}{2}}:
\xi,\zeta\in\E,\psi\in\varpi(\A)\\,\mbox{and}\,\,\psi\bra{\abs{\xi}}=\psi\bra{\abs{\zeta}}=1}.$$
\end{theorem}
\begin{proof} Let $$\varrho=\sup\set{\bra{\abs{\psi\bra{\seq{\zeta,t_1\xi}}}^2+\abs{\psi\bra{\seq{\zeta,t_2\xi}}}^2}^{\frac{1}{2}}:
\xi,\zeta\in\E,\psi\in\varpi(\A)\\,\mbox{and}\,\,\psi\bra{\abs{\xi}}=\psi\bra{\abs{\zeta}}=1}.$$
It is sufficient to prove that $\vertiii{\t}=\varrho$. If $\psi\in\varpi(\A)$ and $\xi, \zeta \in \E$, such that $\psi\bra{\abs{\xi}}=\psi\bra{\abs{\zeta}}=1$,
then by using the Cauchy-Schwartz inequality, we get
\begin{eqnarray*}
  \bra{\abs{\psi\bra{\seq{\zeta,t_1\xi}}}^2+\abs{\psi\bra{\seq{\zeta,t_2\xi}}}^2}^{\frac{1}{2}} &\leq&
\bra{\psi\bra{\seq{t_1\xi,t_1\xi}}\psi\bra{\seq{\zeta,\zeta}}+\psi\bra{\seq{t_2\xi,t_2\xi}}\psi\bra{\seq{\zeta,\zeta}}}^{\frac{1}{2}} \\
   &\leq&\bra{\psi\bra{|t_1\xi|^2}\psi\bra{|\zeta|^2}+\psi\bra{|t_2\xi|^2}\psi\bra{|\zeta|^2}}^{\frac{1}{2}}\\
   &=&
\bra{\psi\bra{|t_1\xi|^2}+\psi\bra{|t_2\xi|^2}}^{\frac{1}{2}}\leq \vertiii{\t}
\end{eqnarray*}
and so $\varrho\leq \vertiii{\t}$.\\
For $\psi\in \varpi(\A)$ and $\xi\in\E$ with $\psi\bra{|\xi|}=1$, we have
$$\psi\bra{|t_j\xi|}^4=\psi\bra{|t_j\xi|^2}=\psi\bra{\seq{t_j\xi,t_j\xi}}=\psi\bra{|t_j\xi|}^2
\psi\bra{\seq{\frac{t_j\xi}{\psi\bra{|t_j\xi|}}},t_j\xi}^2,$$
where assume that $\psi\bra{\abs{t_j\xi}}\neq 0$ for all $j=1,2$. Thus
\begin{eqnarray*}
  \bra{\abs{\psi\bra{t_1\xi}}^2+\abs{\psi\bra{t_2\xi}}^2}^{\frac{1}{2}} &=&\bra{\psi\bra{\seq{\frac{t_1\xi}{\psi\bra{|t_1\xi|}}, t_1\xi}}^2
  +\psi\bra{\seq{\frac{t_2\xi}{\psi\bra{|t_2\xi|}}, t_2\xi}}^2}^{\frac{1}{2}}  \\
   &=& \bra{\abs{\psi\bra{\zeta,t_1\xi}}^2+\abs{\psi\bra{\zeta,t_2\xi}}^2}^{\frac{1}{2}}\leq \varrho
\end{eqnarray*}
and hence $\vertiii{\t}\leq \varrho$.
\end{proof}
Recall that an operator $\t=(t_1,t_2)\in\L^n(\E)$ is said to be self-adjoint if
$\t^*=(t_1^*,t_2^*)=(t_1,t_2)=\t$.
\begin{theorem} If $\t=\bra{t_1,t_2}\in\L^n(\E)$ is self-adjoint, then
$$\vertiii{\t}=\sup\set{\bra{\abs{\psi\bra{\seq{\xi,t_1\xi}}}^2+\abs{\psi\bra{\seq{\xi,t_2\xi}}}^2}^{\frac{1}{2}}
:\xi\in\E,\psi\in\varpi(\A)\,\,\mbox{and}\,\,\psi\bra{|\xi|}=1}.$$
\end{theorem}
\begin{proof} Let $M=\sup\set{\bra{\abs{\psi\bra{\seq{\xi,t_1\xi}}}^2+\bra{\abs{\psi\bra{\seq{\xi,t_j\xi}}}^2}}^{\frac{1}{2}}:
\xi\in\E,\psi\in\varpi(\A)\,\,\mbox{and}\,\,\psi\bra{|\xi|}=1}$. If
$\psi\in\varpi(\A)$ and $\t$ is a self-adjoint, then by using the Cauchy-Schwartz inequality
$$\abs{\psi\bra{\seq{\xi,t_j\xi}}}^2\leq \psi\bra{\seq{t_j\xi,t_j\xi}}\psi\bra{\seq{\xi,\xi}}=\psi\bra{|t_j\xi|^2}\psi\bra{|\xi|^2}$$
for all $j=1,2$. If $\psi\bra{|\xi|}=1$, then
\begin{eqnarray*}
  \bra{\abs{\psi\bra{\seq{\xi,t_1\xi}}}^2+\abs{\psi\bra{\seq{\xi,t_2\xi}}}^2}^{\frac{1}{2}}
  &\leq& \bra{\psi\bra{|t_1\xi|^2}+\psi\bra{|t_2\xi|^2}}^{\frac{1}{2}}\\
   &\leq&\vertiii{\t}.
\end{eqnarray*}
By taking the supremum over all $\xi\in\E$ with $\psi\bra{|\xi|}=1$, we obtain
$$M\leq \vertiii{\t}.$$
For the converse, let  $\xi,\zeta\in\E$ and $\psi\in\varpi(\A)$. Then
$$\psi\bra{\seq{\xi+\zeta,t_j(\xi+\zeta)}}-\psi\bra{\seq{\xi-\zeta,t_j(\xi-\zeta)}}=4\psi\bra{Re\seq{\zeta,t_j\xi}}\,\,\mbox{for $j=1,2$}.$$
Consequently, by Minkowski's inequality we have
\begin{equation*}
  \left.
    \begin{array}{ll}
     \bra{\abs{\psi\bra{Re\seq{\zeta,t_1\xi}}}^2+\abs{\psi\bra{Re\seq{\zeta,t_2\xi}}}^2}^{\frac{1}{2}}\\
=\bra{\frac{1}{16}\abs{\psi\bra{\seq{\xi+\zeta,t_1(\xi+\zeta)}}-\psi\bra{\seq{\xi-\zeta,t_1(\xi-\zeta)}}}^2
+\frac{1}{16}\abs{\psi\bra{\seq{\xi+\zeta,t_2(\xi+\zeta)}}-\psi\bra{\seq{\xi-\zeta,t_2(\xi-\zeta)}}}^2}^{\frac{1}{2}}  \\
   \leq \frac{1}{4}\bra{\sbra{\abs{\psi\bra{\seq{\xi+\zeta,t_1(\xi+\zeta)}}}+\abs{\psi\bra{\seq{\xi-\zeta,t_1(\xi-\zeta)}}}}^2
   +\sbra{\abs{\psi\bra{\seq{\xi+\zeta,t_2(\xi+\zeta)}}}+\abs{\psi\bra{\seq{\xi-\zeta,t_2(\xi-\zeta)}}}}^2}^{\frac{1}{2}}\\
=\frac{1}{4}\bra{\psi\bra{|\xi+\zeta|^4}\abs{\psi\bra{\seq{\frac{\xi+\zeta}{\psi\bra{|\xi+\zeta|}},
\frac{t_1(\xi+\zeta)}{\psi\bra{|\xi+\zeta|}} }}}^2+\psi\bra{|\xi+\zeta|^4}\abs{\psi\bra{\seq{\frac{\xi+\zeta}{\psi\bra{|\xi+\zeta|}},
\frac{t_2(\xi+\zeta)}{\psi\bra{|\xi+\zeta|}} }}}^2}^{\frac{1}{2}} \\
+\frac{1}{4}\bra{\psi\bra{|\xi-\zeta|^4}\abs{\psi\bra{\seq{\frac{\xi-\zeta}{\psi\bra{|\xi-\zeta|}},\frac{t_1(\xi-\zeta)}{\psi\bra{|\xi-\zeta|}} }}}^2+\psi\bra{|\xi-\zeta|^4}\abs{\psi\bra{\seq{\frac{\xi-\zeta}{\psi\bra{|\xi-\zeta|}},\frac{t_2(\xi-\zeta)}{\psi\bra{|\xi-\zeta|}} }}}^2}^{\frac{1}{2}}\\
=\frac{1}{4}\psi\bra{|\xi+\zeta|^2}\bra{\abs{\psi\bra{\seq{\frac{\xi+\zeta}{\psi\bra{|\xi+\zeta|}},
\frac{t_1(\xi+\zeta)}{\psi\bra{|\xi+\zeta|}} }}}^2+\abs{\psi\bra{\seq{\frac{\xi+\zeta}{\psi\bra{|\xi+\zeta|}},
\frac{t_2(\xi+\zeta)}{\psi\bra{|\xi+\zeta|}} }}}^2}^{\frac{1}{2}} \\
+\frac{1}{4}\psi\bra{|\xi-\zeta|^2}\bra{\abs{\psi\bra{\seq{\frac{\xi-\zeta}{\psi\bra{|\xi-\zeta|}},\frac{t_1(\xi-\zeta)}{\psi\bra{|\xi-\zeta|}} }}}^2+\abs{\psi\bra{\seq{\frac{\xi-\zeta}{\psi\bra{|\xi-\zeta|}},\frac{t_2(\xi-\zeta)}{\psi\bra{|\xi-\zeta|}} }}}^2}^{\frac{1}{2}}
    \end{array}
  \right.
\end{equation*}
Since $\psi\bra{\abs{\frac{\xi+\zeta}{\psi\bra{|\xi+\zeta|}}}}=\frac{\psi\bra{|\xi+\zeta|}}{\psi\bra{|\xi+\zeta|}}=1$, we obtain
\begin{equation*}
  \left.
    \begin{array}{ll}
     \bra{\abs{\psi\bra{\seq{\frac{\xi+\zeta}{\psi\bra{|\xi+\zeta|}},
\frac{t_1(\xi+\zeta)}{\psi\bra{|\xi+\zeta|}} }}}^2+\abs{\psi\bra{\seq{\frac{\xi+\zeta}{\psi\bra{|\xi+\zeta|}},
\frac{t_2(\xi+\zeta)}{\psi\bra{|\xi+\zeta|}} }}}^2}^{\frac{1}{2}}\leq M \,\,\mbox{and}\\
\bra{\abs{\psi\bra{\seq{\frac{\xi-\zeta}{\psi\bra{|\xi-\zeta|}},
\frac{t_1(\xi-\zeta)}{\psi\bra{|\xi-\zeta|}} }}}^2+\abs{\psi\bra{\seq{\frac{\xi-\zeta}{\psi\bra{|\xi-\zeta|}},
\frac{t_2(\xi-\zeta)}{\psi\bra{|\xi-\zeta|}} }}}^2}^{\frac{1}{2}}\leq M .
    \end{array}
  \right.
\end{equation*}
Hence
\begin{equation*}
  \left.
    \begin{array}{ll}
      \bra{ \abs{\psi\bra{Re\seq{\zeta,t_1\xi}}}^2+\abs{\psi\bra{Re\seq{\zeta,t_2\xi}}}^2}^{\frac{1}{2}}\leq \frac{1}{4}M\bra{\psi\bra{|\xi+\zeta|^2}+\psi\bra{|\xi-\zeta|^2}}\\
      =\frac{1}{4}M\psi\bra{2|\xi|^2+2|\zeta|^2}=\frac{1}{2}M\psi\bra{|\xi|^2+|\zeta|^2}.
    \end{array}
  \right.
\end{equation*}
If $y=\frac{t_j\xi}{\psi\bra{|t_j\xi|}}$ and $\psi\bra{|\xi|}=1$, then
\begin{eqnarray*}
   \bra{ \abs{\psi\bra{Re\seq{\zeta,t_1\xi}}}^2+\abs{\psi\bra{Re\seq{\zeta,t_2\xi}}}^2}^{\frac{1}{2}}&\leq& \frac{1}{2}M\psi\bra{|\xi|^2+\abs{\frac{t_j\xi}{\psi\bra{|t_j\xi|}}}^2} \\
   &=&\frac{M}{2}\psi\bra{|\xi|^2+\frac{|t_j\xi|^2}{\psi\bra{|t_j\xi|^2}}}\\
&=&\frac{M}{2}\bra{\psi\bra{|\xi|^2}+\frac{\psi\bra{|t_j\xi|^2}}{\psi\bra{|t_j\xi|^2}}}=M.
\end{eqnarray*}
Therefore
\begin{eqnarray*}
 &&\bra{\abs{\frac{1}{\psi\bra{|t_j\xi|}}\psi\bra{Re\seq{t_j\xi,t_1\xi}}}^2
+\abs{\frac{1}{\psi\bra{|t_j\xi|}}\psi\bra{Re\seq{t_2\xi,t_j\xi}}}^2}^{\frac{1}{2}}\\
 &&=
 \bra{ \abs{\frac{1}{\psi\bra{|t_1\xi|}}Re\bra{\psi\bra{|t_1\xi|^2}}}^2
+\abs{\frac{1}{\psi\bra{|t_2\xi|}}Re\bra{\psi\bra{|t_2\xi|^2}}}^2}^{\frac{1}{2}} \\
 &&= \bra{ \abs{\frac{1}{\psi\bra{|t_1\xi|}}\psi\bra{|t_1\xi|^2}}^2+\abs{\frac{1}{\psi\bra{|t_2\xi|}}\psi\bra{|t_2\xi|^2}}^2}^{\frac{1}{2}}\\
&=&\bra{ \abs{\psi\bra{|t_1\xi|}}^2+\abs{\psi\bra{|t_2\xi|}}^2}^{\frac{1}{2}}\leq M
\end{eqnarray*}
and so $\vertiii{\t}\leq M$. The proof of the theorem is complete.
\end{proof}
The following results are very useful in the sequel which can be found in \cite{Mog}.
\begin{lemma}\label{zero1}\label{Euc1} Let $t\in\LE$ and $\psi\in\varpi(\A)$. Then the following are equivalent:
\begin{enumerate}
  \item [(a)] $\psi\bra{\seq{\xi,t\xi}}=0$ for every $\xi\in\E$ with $\psi\bra{|\xi|}=1$;
  \item [(b)] $\psi\bra{\seq{\xi,t\xi}}=0$ for every $\xi\in\E$.
\end{enumerate}
\end{lemma}
\begin{lemma}\label{zero2} Let $t\in\LE$, then $t=0$ if and only if $\psi\bra{\seq{\xi,t\xi}}=0$ for every $\xi\in\E$
and $\psi\in\varpi(\A)$.
\end{lemma}
\begin{lemma} Let $\t=(t_1,t_2)\in\L^n(\E)$. Then for every $\psi\in\varpi(\A)$ and $\xi\in\E$,
  \begin{equation}\label{important99}
    \bra{\abs{\psi\bra{\seq{\xi,t_1\xi}}}^2+\abs{\psi\bra{\seq{\xi,t_2\xi}}}^2}^{\frac{1}{2}}\leq w_e(\t)\psi\bra{\abs{\xi}^2}.
  \end{equation}
\end{lemma}
\begin{proof} For every $\psi\in\varpi(\A)$ and $\xi\in\E$, we have
$$\frac{1}{\psi\bra{|\xi|^4}}\abs{\psi\bra{\seq{\xi,t_j\xi}}}^2=\abs{\psi\bra{\seq{\frac{x}{\psi\bra{|\xi|}},\frac{t_j\xi}{\psi\bra{|\xi|}}}}}^2$$
 for all $j=1,2$. Hence
$$\bra{\frac{1}{\psi\bra{|\xi|^4}}\abs{\psi\bra{\seq{\xi,t_1\xi}}}^2+\frac{1}{\psi\bra{|\xi|^4}}\abs{\psi\bra{\seq{\xi,t_2\xi}}}^2}^{\frac{1}{2}}\leq w_e(\t)$$
and so
$$\bra{\abs{\psi\bra{\seq{\xi,t_1\xi}}}^2+\abs{\psi\bra{\seq{\xi,t_2\xi}}}^2}^{\frac{1}{2}}\leq w_e(\t)\psi\bra{\abs{\xi}^2}.$$
\end{proof}
\begin{theorem}\label{EquivNE} If $\t=(t_1,t_2)\in\L^n(\E)$, then
\begin{equation}\label{equivalent1}
  w_e(\t)\leq \vertiii{\t}\leq 2\sqrt{2}w_e(\t).
\end{equation}
or equivalently,
\begin{equation}\label{equivalent2}
  \frac{1}{2\sqrt{2}}\vertiii{\t}\leq w_e(\t)\leq \vertiii{\t}.
\end{equation}
Here the constants $\frac{1}{2\sqrt{2}}$ and 1 are best possible.
\end{theorem}
\begin{proof} For every $\psi\in\varpi(\A)$ and $\xi\in\E$ such that $\psi\bra{\abs{\xi}}=1$, by Theorem \ref{extension}, we have
$$\bra{\abs{\psi\bra{\xi,t_1\xi}}^2+\abs{\psi\bra{\xi,t_2\xi}}^2}^{\frac{1}{2}}\leq \vertiii{\t},$$
by taking the supremum, we obtain
$$w_e(\t)\leq \vertiii{\t}.$$
Fix $\xi,\zeta\in\E$ and $\psi\in\varpi(\A)$, we have for all $j=1,2$ that
\begin{eqnarray*}
  4\abs{\psi\bra{\seq{\zeta,t_j\xi}}} &=& |\psi[\seq{\xi+\zeta,t_j(\xi+\zeta)}-\seq{\xi-\zeta,t_j(\xi-\zeta)}\\
&+&i\seq{\xi+i\zeta,t_j(\xi+i\zeta)}-i\seq{\xi-i\zeta,t_j(\xi-i\zeta)}] |\\
   &\leq&\abs{\psi\bra{\seq{\xi+\zeta,t_j(\xi+\zeta)}}}+\abs{\psi\bra{\seq{\xi-\zeta,t_j(\xi-\zeta)}}}\\
&+&\abs{\psi\bra{\seq{\xi+i\zeta,t_j(\xi+i\zeta)}}}
+\abs{\psi\bra{\seq{\xi-i\zeta,t_j(\xi-i\zeta)}}}.
\end{eqnarray*}
Hence for all $j=1,2$, we have
\begin{eqnarray*}
&&4\abs{\psi\bra{\seq{\zeta,t_j\xi}}} \leq
\bra{\abs{\psi\bra{\seq{\xi+\zeta,t_1(\xi+\zeta)}}}^2+\abs{\psi\bra{\seq{\xi+\zeta,t_j(\xi+\zeta)}}}^2}^{\frac{1}{2}}\\
&&+\bra{\abs{\psi\bra{\seq{\xi-\zeta,t_1(\xi-\zeta)}}}^2+\abs{\psi\bra{\seq{\xi-\zeta,t_2(\xi-\zeta)}}}^2}^{\frac{1}{2}}\\
&&+\bra{\abs{\psi\bra{\seq{\xi+i\zeta,t_1(\xi+i\zeta)}}}^2+\abs{\psi\bra{\seq{\xi+i\zeta,t_2(\xi+i\zeta)}}}^2}^{\frac{1}{2}}\\
&&+\bra{\abs{\psi\bra{\seq{\xi-i\zeta,t_1(\xi-i\zeta)}}}^2+\abs{\psi\bra{\seq{\xi-i\zeta,t_2(\xi-i\zeta)}}}^2}^{\frac{1}{2}}\\
&\leq& w_e(\t)\psi\bra{|\xi+\zeta|^2}+w_e(\t)\psi\bra{|\xi-\zeta|^2}+w_e(\t)\psi\bra{|\xi+i\zeta|^2}+w_e(\t)\psi\bra{|\xi-i\zeta|^2}\\
&=&w_e(\t)\sbra{\psi\bra{|\xi+\zeta|^2}+\psi\bra{|\xi-\zeta|^2}+\psi\bra{|\xi+i\zeta|^2}+\psi\bra{|\xi-i\zeta|^2}}\\
&=&w_e(\t)\psi\bra{2|\xi|^2+2|\zeta|^2+2|\xi|^2+2|i\zeta|^2}=4w_e(\t)\bra{\psi(|\xi|^2)+\psi(|\zeta|^2)}.
\end{eqnarray*}
If $\psi(|\xi|)=\psi(|\zeta|)=1$, then
$$\abs{\psi\bra{\seq{\zeta,t_j\xi}}}\leq 2w_e(\t).$$
Now
$$\bra{\abs{\psi\bra{\seq{\zeta,t_1\xi}}}^2+\abs{\psi\bra{\seq{\zeta,t_2\xi}}}^2}^{\frac{1}{2}}\leq \bra{8w_e^2(\t)}^{\frac{1}{2}}=2\sqrt{2}w_e(\t).$$
Taking the supremum over all $\xi,\zeta\in\E$ and $\psi\in\varpi(\A)$ such that $\psi(|\xi|)=\psi(|\zeta|)=1$,we get
$$\vertiii{\t}\leq 2\sqrt{2}w_e(\t).$$
\end{proof}
\begin{remark} (i) Observe that if $\t=(\x,\y)$ is a self adjoint, then it follows from (\ref{equivalent2}) that
\begin{equation}\label{MUNA1}
  \frac{1}{2\sqrt{2}}\vertiii{\x^2+\y^2}^{\frac{1}{2}}\leq w_e(\x,\y)\leq \vertiii{\x^2+\y^2}^{\frac{1}{2}}.
\end{equation}
(ii) If $\t=\x+i\y$ is the cartesian decomposition of $\t$, then
\begin{eqnarray}\label{MUNA2}
  w_e^2(\x,\y)&=&\sup_{\psi\bra{|x|}=1}\sbra{\abs{\psi\bra{\seq{\xi,\x \xi}}}^2+\abs{\psi\bra{\seq{\xi,\y \xi}}}^2}\nonumber \\
   &=&\sup_{\psi\bra{|x|}=1}\abs{\psi\bra{\seq{\xi,\t \xi}}}^2=w^2(\t).
\end{eqnarray}
(iii) If $\t=\x+i\y$ is the cartesian decomposition of $\t$, then
$$\t^*\t+\t\t^*=2(\x^2+\y^2)$$
and hence it follows from (\ref{MUNA1}) that
$$\frac{1}{16}\vertiii{\t^*\t+\t\t^*}\leq w^2(\t)\leq \frac{1}{2}\vertiii{\t^*\t+\t\t^*}.$$
\end{remark}
\begin{theorem} Let $\t=\x+i\y$ be the cartesian decomposition of $\t\in\L(\E)$. Then for every
$\mu,\nu\in\Real$,
\begin{equation}\label{Cart1}
  w_e(\x,\y)=\sup_{\mu^2+\nu^2=1}\vertiii{\mu \x+\nu \y}.
\end{equation}
In particular,
\begin{equation}\label{Cart2}
  \frac{1}{2}\vertiii{\t+\t^*}\leq w_e(\x,\y)\,\,\, \mbox{and}\,\,\,\frac{1}{2}\vertiii{\t-\t^*}\leq w_e(\x,\y).
\end{equation}
\end{theorem}
\begin{proof} First of all, we note that
\begin{equation}\label{Cart3}
  w_e(\x,\y)=\sup_{\theta\in\Real}\vertiii{Re\bra{e^{i\theta}\t}}.
\end{equation}
In fact, $\sup_{\theta\in\Real}\bra{e^{i\theta}\psi\bra{\xi,\t\xi}}^2=\abs{\psi\bra{\seq{\xi,\x\xi}}}^2+\abs{\psi\bra{\seq{\xi,\y\xi}}}^2$ yields to
$$\sup_{\theta\in\Real}\vertiii{Re\bra{e^{i\theta}\t}}=\sup_{\theta\in\Real}w\bra{Re(e^{i\theta}\t)}=w_e(\x,\y).$$
On the other hand, let $\t=\x+i\y$ be the Cartesian decomposition of $\t$. Then
\begin{eqnarray}\label{Cart4}
  Re\bra{e^{i\theta}\t} &=& \frac{e^{i\theta}\t+e^{-i\theta}\t^*}{2}\nonumber \\
   &=&\frac{1}{2}\sbra{\bra{\cos\theta+i\sin\theta}\t+\bra{\cos\theta-i\sin\theta}\t^*}\nonumber\\
&=&\bra{\cos\theta}\bra{\frac{\t+\t^*}{2}}-\bra{\sin\theta}\bra{\frac{\t-\t^*}{2i}}\nonumber\\
&=&\bra{\cos\theta} \x-\bra{\sin\theta}\y.
\end{eqnarray}
Therefore, by putting $\mu=\cos\theta$ and $\nu=\sin\theta$ in (\ref{Cart4}), we obtain (\ref{Cart1}). Especially, by
setting $(\mu,\nu)=(1,0)$ and $(\mu,\nu)=(0,1)$, we reach (\ref{Cart2}).
\end{proof}
\begin{remark} By using (\ref{Cart2}), we get some known inequalities:
\begin{enumerate}
  \item [(a)] $\vertiii{\t}=\vertiii{\x+i\y}\leq \vertiii{\x}+\vertiii{\y}\leq 2w_e(\x,\y)$. Hence we have
$\frac{1}{2}\vertiii{\t}\leq w_e(\x,\y)$.
  \item [(b)] If $\t=\t^*$, then $\t=\x$. Hence we have $\vertiii{\t}=\vertiii{\x}\leq w_e(\x,\y)\leq \vertiii{\t}$
and $w_e(\x,\y)=\vertiii{\t}$.
  \item [(c)] By an easy calculation, we have $\frac{\t^*\t+\t\t^*}{2}=\x^2+\y^2$. Hence
$$\frac{1}{4}\vertiii{\t^*\t+\t\t^*}=\frac{1}{2}\vertiii{\x^2+\y^2}\leq \frac{1}{2}\bra{\vertiii{\x}^2+\vertiii{\y}^2}\leq w_e^2(\x,\y).$$
\item[(d)] Let $\mu,\nu\in\Real$ such that $\mu^2+\nu^2=1$. Then for every $\xi\in\E$ and $\psi\in\varpi(\A)$ such that
$\psi\bra{|\xi|}=1$, we have
\begin{eqnarray*}
  \vertiii{(\mu\x+\nu\y)\xi} &=& \vertiii{\begin{bmatrix} \x & \y \\ 0 & 0 \\ \end{bmatrix}\begin{bmatrix} \mu\xi \\\nu\xi \\ \end{bmatrix}} \leq \vertiii{\begin{bmatrix} \x & \y \\ 0 & 0 \\ \end{bmatrix}}\leq \vertiii{\begin{bmatrix} \x & \y \\ 0 & 0 \\ \end{bmatrix}\begin{bmatrix} \x & 0 \\ \y & 0 \\ \end{bmatrix}}^{\frac{1}{2}}\\
&=&\vertiii{\x^2+\y^2}^{\frac{1}{2}}=\frac{1}{\sqrt{2}}\vertiii{\t^*\t+\t\t^*}^{\frac{1}{2}}.
\end{eqnarray*}
Hence we have
$$w_e^2(\x,\y)=\sup_{\mu^2+\nu^2=1}\vertiii{\mu\x+\nu\y}^2\leq \frac{1}{2}\vertiii{\t^*\t+\t\t^*}.$$
\end{enumerate}
\end{remark}
\begin{theorem}\label{normpf} $w_e:\L^2(\E)\to [0,\infty)$ is  defines a norm which is equivalent to the
norm on $\L^2(\E)$.
\end{theorem}
\begin{proof} Let $w_e(t_1,t_2)=0$. Then for every $\xi\in \E$ and $\psi\in \varpi(\A)$ with $\psi\bra{|\xi|}=1$,
we have
$$\bra{\abs{\psi\bra{\seq{\xi,t_1\xi}}}^2+\abs{\psi\bra{\seq{\xi,t_2\xi}}}^2}^{\frac{1}{2}}=0.$$
Hence $\psi\bra{\seq{\xi,t_k\xi}}=0$  for every $\xi\in\E$ and $\psi\in \varpi(\A)$ with
$\psi\bra{|\xi|}=1$ and $k=1,2$. By Lemma \ref{zero1}, $\psi\bra{\seq{\xi,t_k\xi}}=0$  for every $\xi\in\E$ and $\psi\in \varpi(\A)$ with
$\psi\bra{|\xi|}=1$, and by Lemma \ref{zero2}, $t_k=0$ for $k=1,2$. For every $\mu\in\c,$
\begin{eqnarray*}
  w_e(\mu t_1,\mu t_2)&=&\sup_{\psi\bra{|\xi|}=1}\bra{\abs{\psi\bra{\seq{x,\mu t_1x}}}^2+\abs{\psi\bra{\seq{x,\mu t_2x}}}^2}^{\frac{1}{2}} \\
  &=&\abs{\mu}\bra{\abs{\psi\bra{\seq{\xi,t_1\xi}}}^2+\abs{\psi\bra{\seq{\xi,t_2\xi}}}^2}^{\frac{1}{2}}=\abs{\mu}w_e(t_1,t_2).
\end{eqnarray*}
Let $t_1,t_2, s_1,s_2\in\LE$. For every $\xi\in\E$ and $\psi\in \varpi(\A)$ with $\psi\bra{|\xi|}=1$, and by Minkowski's inequality  we have
\begin{eqnarray*}
  &&\sbra{\abs{\psi\bra{\seq{\xi,(t_1+s_1)\xi}}}^2+\abs{\psi\bra{\seq{\xi,(t_2+s_2)\xi}}}^2}^{\frac{1}{2}}\\
&&\leq \sbra{\bra{\abs{\psi\bra{\seq{\xi,t_1\xi}}}+\abs{\psi\bra{\seq{x,s_1x}}}}^2+
\bra{\abs{\psi\bra{\seq{\xi,t_2\xi}}}+\abs{\psi\bra{\seq{x,s_2x}}}}^2}^{\frac{1}{2}}\\
&\leq& \sbra{\abs{\psi\bra{\seq{\xi,t_1\xi}}}^2+\abs{\psi\bra{\seq{\xi,t_2\xi}}}^2}^{\frac{1}{2}}+
\sbra{\abs{\psi\bra{\seq{\xi,s_1\xi}}}^2+\abs{\psi\bra{\seq{\xi,s_2\xi}}}^2}^{\frac{1}{2}}
\end{eqnarray*}
By taking supremum over all $\psi\bra{|\xi|}=1$,
$$w_e\bra{t_1+s_1,t_2+s_2}\leq w_e(t_1,t_2)+w_e(s_1,s_2).$$
Finally, the equivalence with the norm on $\L^2(\E)$ follows from Theorem \ref{EquivNE}.
\end{proof}
\section{Main results}
\begin{theorem}\label{Karakash1} Let $t,s\in\LE$. Then
\begin{equation}\label{Prof1}
  \frac{1}{\sqrt{2}}\sbra{w(t^2+s^2)}^{\frac{1}{2}}\leq w_e(t,s)\leq \vertiii{t^*t+s^*s}^{\frac{1}{2}}.
\end{equation}
\end{theorem}
\begin{proof} We follow a similar argument to the one from \cite{Mog}. For every $x\in\E$ with $\psi\bra{|x|}=1$ and $\psi\in\varpi(\A)$, we have
\begin{eqnarray}\label{Prof2}
 \abs{\psi\bra{\seq{x,tx}}}^2+\abs{\psi\bra{\seq{x,sx}}}^2&\geq& \frac{1}{2}\bra{\abs{\psi\bra{\seq{x,tx}}+\abs{\psi\bra{\seq{x,tx}}}}}^2\nonumber\\
  &\geq& \frac{1}{2}\abs{\psi\bra{\seq{x,(t\pm s)x}}}^2.
\end{eqnarray}
Taking the supremum in (\ref{Prof2}), we deduce
  \begin{equation}\label{prof3}
    w_e^2(t,s)\geq \frac{1}{2}w^2(t\pm s).
  \end{equation}
Utilising the inequality (\ref{prof3}) and the properties of the numerical radius, we have successively:
\begin{eqnarray*}
  2w_e^2(t,s)&\geq & \frac{1}{2}\sbra{w^2(t+s)+w^2(t-s)}\\
   &\geq& \frac{1}{2}\sbra{w\bra{(t+s)^2}+w\bra{(t-s)^2}}\\
&\geq&  \frac{1}{2}w\bra{(t+s)^2+(t-s)^2}\\
&=&w(t^2+s^2).
\end{eqnarray*}
which gives the desired inequality (\ref{Prof1}).
\end{proof}
\begin{corollary} For any two self-adjoint bounded linear operators $t,s$ on $\E$, we have
\begin{equation}\label{prof4}
  \frac{1}{\sqrt{2}}\vertiii{t^2+s^2}^{\frac{1}{2}}\leq w_e(t,s)\leq \vertiii{t^2+s^2}^{\frac{1}{2}}.
\end{equation}
\end{corollary}
\begin{remark} Observe that for the case that $t,s$ are the self-adjoint
operators in the Cartesian decomposition of $a$, exactly the lower bound obtained by  Moghaddam and  Mirmostafaee
in \cite[Theorem 3.2]{Mog} for the numerical radius $w(a)$. Moreover, since $\frac{1}{4}$
 is a sharp constant in Moghaddam and  Mirmostafaee inequality, it follows that $\frac{1}{\sqrt{2}}$ is also the best possible constant in (\ref{prof4}) and (\ref{Prof1}), respectively.
\end{remark}
\begin{corollary} For any $a\in\LE$ and $\mu,\nu\in\c$, we have
\begin{equation}\label{prof5}
  \frac{1}{2}w\bra{\mu^2 a^2+\nu^2(a^*)^2}\leq \bra{|\mu|^2+|\nu|^2}w^2(a)\leq \vertiii{|\mu|^2a^*a+|\nu|^2aa^*}.
\end{equation}
\end{corollary}
\begin{proof} In Theorem \ref{Karakash1}, put $t=\mu a$ and $s=\nu a^*$, we get
$$w_e^2(t,s)=\bra{|\mu|^2+|\nu|^2}w^2(a)$$
and
$$w^2(t^2+s^2)=w\bra{|\mu|^2a^*a+|\nu|^2aa^*},$$
which, by (\ref{Prof1}) implies the desired result (\ref{prof5}).
\end{proof}
\begin{remark}(i) If we choose in (\ref{prof5}) $\mu=\nu\neq 0$, then we obtain the following inequality
\begin{equation}\label{prof6}
  \frac{1}{4}\vertiii{a^2+(a^*)^2}\leq w^2(a)\leq \frac{1}{2}\vertiii{a^*a+aa^*}
\end{equation}
for any operator $a\in\LE$.\\
(ii) If we choose in (\ref{prof5}) $\mu=1,\nu=i$, then we obtain the following inequality
\begin{equation}\label{prof7}
  \frac{1}{4}\vertiii{a^2-(a^*)^2}\leq w^2(a)
\end{equation}
for any operator $a\in\LE$.
\end{remark}
\begin{theorem} For any $t,s\in\LE$, we have
\begin{equation}\label{prof8}
  \frac{1}{\sqrt{2}}\max\{w(t+s),w(t-s)\}\leq w_e(t,s)\leq \frac{1}{\sqrt{2}}\sbra{w^2(t+s)+w^2(t-s)}^{\frac{1}{2}}.
\end{equation}
The constant $\frac{1}{\sqrt{2}}$ is sharp in both inequalities.
\end{theorem}
\begin{proof} The first inequality follows from (\ref{prof3}).
For the second inequality, we observe that
\begin{equation}\label{prof9}
  \abs{\psi\bra{\seq{x,tx}}\pm \psi\bra{\seq{x,sx}}}^2\leq w^2(t\pm s)
\end{equation}
for every $x\in\E$ with $\psi\bra{|x|}=1$ and $\psi\in \varpi(\A)$.\\
The inequality (\ref{prof9}) and the parallelogram identity for complex numbers give
\begin{eqnarray}\label{prof10}
  2\sbra{\abs{\psi\bra{\seq{x,tx}}}^2+ \abs{\psi\bra{\seq{x,sx}}}^2}&=& \abs{\psi\bra{\seq{x,tx}}-\psi\bra{\seq{x,sx}}}^2
+\abs{\psi\bra{\seq{x,tx}}+\psi\bra{\seq{x,sx}}}^2\nonumber \\
   &\leq&w^2(t+s)+w^2(t-s)
\end{eqnarray}
for every $x\in\E$ with $\psi\bra{|x|}=1$ and $\psi\in \varpi(\A)$.\\
Taking the supremum in (\ref{prof9}) over all vector $x\in\E$ with $\psi\bra{|x|}=1$, we deduce the desired result (\ref{prof8}).\\
The fact that $\frac{1}{\sqrt{2}}$ is the best possible constant follows from the fact that for $t=s\neq 0$ one
would obtain the same quantity $\sqrt{2}w(t)$ in all terms of (\ref{prof8}).
\end{proof}
\begin{corollary} For any two self-adjoint operators $t,s\in\LE$, we have
\begin{equation}\label{prof8}
  \frac{1}{\sqrt{2}}\max\{\vertiii{t+s},\vertiii{t-s}\}\leq w_e(t,s)\leq \frac{1}{\sqrt{2}}\sbra{\vertiii{t+s}^2+\vertiii{t-s}^2}^{\frac{1}{2}}.
\end{equation}
The constant $\frac{1}{\sqrt{2}}$ is sharp in both inequalities.
\end{corollary}
\begin{theorem}\label{prop1} Let $t,s\in\LE$. Then
\begin{equation}\label{prof16}
  w_e(t,s)\leq \sbra{w^2(t-s)+2w(t)w(s)}^{\frac{1}{2}}.
\end{equation}
\end{theorem}
\begin{proof} For any $x\in\E$ with $\psi\bra{|x|}=1$ and $\psi\in\varpi(\A)$, we have
$$\abs{\psi\bra{\seq{x,tx}}}^2-2Re\sbra{\psi\bra{\seq{x,tx}}\overline{\psi\bra{\seq{x,sx}}}}+
\abs{\psi\bra{\seq{x,sx}}}^2=\abs{\psi\bra{\seq{x,tx}}-\psi\bra{\seq{x,sx}}}^2\leq w^2(t-s),$$
giving
  \begin{eqnarray}\label{prof17}
    \abs{\psi\bra{\seq{x,tx}}}^2+\abs{\psi\bra{\seq{x,sx}}}^2 &\leq& w^2(t-s)+ 2Re\sbra{\psi\bra{\seq{x,tx}}\overline{\psi\bra{\seq{x,sx}}}}\nonumber\\
     &\leq&w^2(t-s)+2\abs{\psi\bra{\seq{x,tx}}}\abs{\psi\bra{\seq{x,sx}}}
  \end{eqnarray}
for any $x\in\E$ with $\psi\bra{|x|}=1$ and $\psi\in\varpi(\A)$.\\
Taking the supremum in (\ref{prof17}) over all $x\in\E$ with $\psi\bra{|x|}=1$ and $\psi\in\varpi(\A)$, we deduce the required inequality (\ref{prof16}).
\end{proof}
In particular, if $t$ and $s$ are self-adjoint operators, then
\begin{equation}\label{prof18}
  w_e(t,s)\leq \sbra{\vertiii{t-s}^2+\vertiii{t+s}^2}^{\frac{1}{2}}.
\end{equation}
The following result provides a different upper bound for the Euclidean operator radius than
(\ref{prof16}).
\begin{theorem}\label{prop2} Let $t,s\in\LE$. Then we have
\begin{equation}\label{prof19}
  w_e(t,s)\leq \sbra{2\min\{w^2(t),w^2(s)\}+w(t-s)w(t+s)}^{\frac{1}{2}}.
\end{equation}
\end{theorem}
\begin{proof} Utilising the parallelogram identity (\ref{prof10}), we have, by taking the supremum over $x\in\E$,
$\psi\bra{|x|}=1$ and $\psi\in\varpi(\A)$, that
\begin{equation}\label{prof20}
  2w_e^2(t,s)=w_e^2(t-s,t+s).
\end{equation}
Now, if we apply Theorem \ref{prop1} for $t-s,t+s$ instead of $t$ and $s$, then we can state
$$w_e^2(t-s,t+s)\leq 4w^2(s)+2w(t-s)w(t+s)$$
giving
\begin{equation}\label{prof21}
  w_e^2(t,s)\leq 2w^2(s)+w(t-s)w(t+s).
\end{equation}
Now, if in (\ref{prof21}) we swap the $s$ with $t$ then we also have
\begin{equation}\label{prof22}
  w_e^2(t,s)\leq 2w^2(t)+w(t-s)w(t+s).
\end{equation}
The conclusion follows now by (\ref{prof21}) and (\ref{prof22}).
\end{proof}
A different upper bound for the Euclidean operator radius is incorporated in the following.
\begin{theorem}\label{Rahma99} Let $t,s\in\LE$. Then we have
\begin{equation}\label{prof23}
  w_e^2(t,s)\leq \max\{\vertiii{t}^2,\vertiii{s}^2\}+w(s^*t).
\end{equation}
\end{theorem}
\begin{proof}Firstly, let us observe that for any $y, u, v \in\E$ and $\psi\in\varpi(\A)$ we have successively
\begin{eqnarray}\label{prof24}
 && \vertiii{\psi\bra{\seq{y,u}}u+\psi\bra{\seq{y,v}}v}^2 =\nonumber\\
&& \abs{\bra{\seq{y,u}}}^2\psi\bra{|u|^2}+ \abs{\bra{\seq{y,v}}}^2
\psi\bra{|v|^2}+2Re\sbra{\psi\bra{\seq{y,u}}\overline{\psi\bra{\seq{y,v}}}\psi\bra{\seq{u,v}}}\nonumber\\
 &&  \leq \abs{\bra{\seq{y,u}}}^2\psi\bra{|u|^2}+ \abs{\bra{\seq{y,v}}}^2
\psi\bra{|v|^2}+2\abs{\psi\bra{\seq{y,u}}}\abs{\psi\bra{\seq{y,v}}}\abs{\psi\bra{\seq{u,v}}}\nonumber\\
&&\leq \abs{\bra{\seq{y,u}}}^2\psi\bra{|u|^2}+ \abs{\bra{\seq{y,v}}}^2\psi\bra{|v|^2}+\sbra{\abs{\psi\bra{\seq{y,u}}}^2+\abs{\psi\bra{\seq{y,v}}}^2}
\abs{\psi\bra{\seq{u,v}}}\nonumber\\
&&\leq \sbra{\abs{\psi\bra{\seq{y,u}}}^2+\abs{\psi\bra{\seq{y,v}}}^2}\bra{\max\{\psi\bra{|u|^2},\psi\bra{|v|^2}\}+
\abs{\psi\bra{\seq{u,v}}}}.
\end{eqnarray}
On the other hand,
\begin{eqnarray}\label{prof25}
  &&\sbra{\abs{\psi\bra{\seq{y,u}}}^2+\abs{\psi\bra{\seq{y,v}}}^2}^2=\sbra{\psi\bra{\seq{y,u}}\psi\bra{\seq{u,y}}
+\psi\bra{\seq{y,v}}\psi\bra{\seq{v,y}}}^2\nonumber \\
&&=\sbra{\psi\bra{\seq{y,\psi\bra{\seq{y,u}}u+\psi\bra{\seq{y,v}}v}}}^2\nonumber\\
&&\leq \psi\bra{|y|^2}\psi\bra{\abs{\psi\bra{\seq{y,u}}u+\psi\bra{\seq{y,v}}v}^2}.
\end{eqnarray}
for any $u,v,y\in\E$.\\
Making use of (\ref{prof24}) and (\ref{prof25}) we deduce that
\begin{equation}\label{prof26}
  \abs{\psi\bra{\seq{y,u}}}^2+\abs{\psi\bra{\seq{y,v}}}^2\leq \psi\bra{|y|^2}\bra{\max\{\psi\bra{|u|^2},\psi\bra{|v|^2}\}+
\abs{\psi\bra{\seq{u,v}}}}
\end{equation}
for any $u,v,y\in\E$.\\
Now, if we apply the inequality (\ref{prof26}) for $y =x$, $u=tx$, $v=sx$, $x\in\E$, $\psi\bra{|x|}=1$, then we
can state that
\begin{equation}\label{prof27}
  \abs{\psi\bra{\seq{x,tx}}}^2+\abs{\psi\bra{\seq{x,sx}}}^2\leq \psi\bra{|x|^2}\bra{\max\{\vertiii{tx}^2,\vertiii{sx}^2\}+
\abs{\psi\bra{\seq{tx,sx}}}}
\end{equation}
for any $x\in\E$ with $\psi\bra{|x|}=1$ and $\psi\in\varpi(\A)$.\\
Taking the supremum over $x\in\E$ with $\psi\bra{|x|}=1$ and $\psi\in\varpi(\A)$, we deduce the desired result (\ref{prof23}).
\end{proof}
\begin{remark} In Theorem \ref{Rahma99},  if $t$ and $s$ are self-adjoint and $t=s$,  then
the inequality (\ref{prof23}) becomes
$$w_e(t,t)\leq \sqrt{2}\vertiii{t}.$$
Hence the inequality (\ref{prof23}) is sharp.
\end{remark}
If information about the sum and the difference of the operators B and C are available, then
one may use the following result.
\begin{corollary} Let $t,s\in\LE$. Then
\begin{equation}\label{prof28}
  w_e^2(t,s)\leq \frac{1}{2}\sbra{\max\{\vertiii{t-s}^2,\vertiii{t+s}^2\}+w[(t+s)(t^*-s^*)]}.
\end{equation}
  The constant $\frac{1}{2}$ is best possible.
\end{corollary}
\begin{proof} Follows by the inequality (\ref{prof23}) written for $t+s$ and $t-s$ instead of $t$ and $s$ and by
utilising the identity (\ref{prof20}).
The fact that $\frac{1}{2}$ is best possible in (\ref{prof28}) follows by the fact that for $s=t$, $t$ a self-adjoint
operator, we get in both sides of the inequality (\ref{prof28}) the quantity $2\vertiii{t}^2$.
\end{proof}
\begin{corollary} Let  $a\in\LE$. Then
\begin{equation}\label{prof29}
  w^2(t)\leq \frac{1}{4}\sbra{\max\{\vertiii{t-t^*}^2,\vertiii{t+t^*}^2\}+w[(t-t^*)(t+t^*)]}.
\end{equation}
The constant $\frac{1}{4}$ is best possible.
\end{corollary}
\begin{proof} If $t=\frac{a+a^*}{2}$ and $s=\frac{a-a^*}{2i}$ are the cartesian decomposition of $a\in\LE$,
then
$$w_e^2(t,s)=w^2(a)$$
and
$$w(s^*t)=\frac{1}{4}\sbra{(a^*+a)(a^*-a)}.$$
 Utilizing (\ref{prof23}), we deduce (\ref{prof29}).
\end{proof}
\begin{remark} If we choose in (\ref{prof23}), $t=a$ and $s=a^*$, $a\in\LE$, then we can state that
\begin{equation}\label{prof30}
  w^2(a)\leq \frac{1}{2}\sbra{\vertiii{a}^2+w(a^2)}.
\end{equation}
The constant $\frac{1}{2}$ is best possible.
\end{remark}
the following upper bound for the Euclidean radius involving different composite
operators
\begin{theorem} Let $t,s\in\LE$. Then
\begin{equation}\label{prof31}
   w_e^2(t,s)\leq \frac{1}{2}\sbra{\vertiii{t^*t+s^*s}+\vertiii{t^*t-s^*s}}+w(s^*t).
\end{equation}
\end{theorem}
\begin{proof} We use (\ref{prof27}) to write that
\begin{equation}\label{prof32}
  \abs{\psi\bra{\seq{x,tx}}}^2+\abs{\psi\bra{\seq{x,sx}}}^2\leq \frac{1}{2}\sbra{\vertiii{tx}^2+\vertiii{sx}^2
+\abs{\vertiii{tx}^2-\vertiii{sx}^2}}+\abs{\psi\bra{\seq{sx,tx}}}
\end{equation}
for any $x\in\E$ with $\psi\bra{|x|}=1$ and $\psi\in\varpi(\A)$.\\
Since $\vertiii{tx}^2=\psi\bra{\seq{x,t^*tx}}$ and $\vertiii{sx}^2=\psi\bra{\seq{x,s^*sx}}$, then (\ref{prof32}) can be written
as
\begin{eqnarray}\label{prof33}
  \abs{\psi\bra{\seq{x,tx}}}^2+\abs{\psi\bra{\seq{x,sx}}}^2&\leq&\frac{1}{2}\sbra{\psi\bra{\seq{x,(t^*t+s^*s)x}}
+\abs{\psi\bra{\seq{x,(t^*t-s^*s)x}}}}\nonumber \\
   &+& \abs{\psi\bra{\seq{sx,tx}}}
\end{eqnarray}
for any $x\in\E$ with $\psi\bra{|x|}=1$ and $\psi\in\varpi(\A)$.\\
Taking the supremum in (\ref{prof33}) over $x\in\E$ with $\psi\bra{|x|}=1$ and $\psi\in\varpi(\A)$. and noticing that the operators
$t^*t\pm s^*s$ are self-adjoint, we deduce the desired result (\ref{prof31}).
\end{proof}
\begin{corollary} for any operators $t,s\in\LE$, we have
\begin{equation}\label{prof34}
  w_e^2(t,s)\leq \frac{1}{2}\sbra{\vertiii{t^*t+s^*s}+\vertiii{t^*s+s^*t}+w[(t+s)(t^*-s^*)]}.
\end{equation}
The constant $\frac{1}{2}$ is best possible.
\end{corollary}
\begin{proof} If we write (\ref{prof31}) for $t+s,t-s$ instead of $t,s$ and perform the required calculations
then we get
$$w_e^2(t+s,t-s)\leq \frac{1}{2}\sbra{2\vertiii{t^*t+s^*s}+2\vertiii{t^*s+s^*t}}+w[(t+s)(t^*-s^*)],$$
which, by the identity (\ref{prof20}) is clearly equivalent with (\ref{prof34}).\\
Now, if we choose in (\ref{prof34}) $t=s$, then we get the inequality $w(t)\leq \vertiii{t}$, which is a sharp
inequality.
\end{proof}
\begin{corollary} If $t,s\in\LE$ are self-adjoint, then
\begin{equation}\label{prof35}
  w_e^2(t,s)\leq \frac{1}{2}\sbra{\vertiii{t^2+s^2}+\vertiii{t^2-s^2}}+w(st).
\end{equation}
\end{corollary}
\begin{remark}(i) We observe that, if $t$ and $s$ are chosen to be the Cartesian decomposition for the bounded
linear operator $a$, then we can get from (\ref{prof35}) that
\begin{equation}\label{prof36}
  w^2(a)\leq \frac{1}{4}\sbra{\vertiii{t^*t+tt^*}+\vertiii{t^2+(t^*)^2}+w\bra{(a+a^*)(a^*-a)}}.
\end{equation}
The constant $\frac{1}{4}$ is best possible. This follows by the fact that for $a$ a self-adjoint operator, we
obtain in both sides of (\ref{prof36}) the same quantity $\vertiii{a}^2$.\\
(ii) Now, if we choose in (\ref{prof31}) $t-a$ and $s=a^*$, $a\in\LE$, then we get
\begin{equation}\label{prof37}
  w^2(a)\leq \frac{1}{4}\sbra{\vertiii{a^*a+a^*a}+\vertiii{a^*a-aa^*}}+\frac{1}{2}w(a^2).
\end{equation}
This inequality is sharp. The equality holds if, for instance, we  assume that $a$ is normal, i.e., $a^*a=aa^*$. In this case we get in both sides of (\ref{prof37}) the quantity $\vertiii{a}^2$, since for normal operators, $w(a^2)=w^2(a)=\vertiii{a}^2$.
\end{remark}
{\bf Author Contributions:}   The author have read and agreed to the published version of the manuscript.\\
{\bf Funding:} No funding is applicable\\
{\bf Institutional Review Board Statement:} Not applicable.\\
{\bf Informed Consent Statement:} Not applicable.
{\bf Data Availability Statement:} Not applicable.\\
{\bf Conflicts of Interest:} The authors declare no conflict of interest.

\bibliographystyle{unsrtnat}
\bibliography{references}  






\end{document}